\documentclass[10pt]{article}
\usepackage{float}
\usepackage{graphicx}
\usepackage{epstopdf}
\usepackage{cite}
\usepackage{mathrsfs}
\usepackage{amsfonts}
\usepackage{amsmath}
\usepackage{amsfonts,amssymb,color}
\usepackage{dsfont}
\usepackage{curves}
\usepackage{mathrsfs}
\usepackage{pifont}
\usepackage{amssymb}
\usepackage{tikz}
\usepackage{hyperref}
\usepackage{enumitem}
\usepackage{CJK}
\usepackage[english]{babel}
\usepackage{times}
\usepackage[T1]{fontenc}
\usepackage{indentfirst}
\usepackage{subfigure}

\usepackage{subfigure}
\usepackage[T1]{fontenc}
\usepackage[utf8]{inputenc}
\usepackage{authblk}
\usepackage[misc]{ifsym}
\newtheorem{theorem}{\bf Theorem}[section]

\newtheorem{corollary}[theorem]{\bf Corollary}

\newtheorem{problem}{\bf Problem}[section]

\newcommand{\qed}{\hfill\rule{0.5em}{0.809em}}

\def\pf{\noindent {\it Proof. }}

\def\qed{\hfill \rule{4pt}{7pt}}

\date{}
\begin{document}


\title{ On the existence of tripartite graphs and $n$-partite graphs
 \thanks{Supported by National Natural
 	Science Foundation of China (Grant No. 11961019) and Hainan Provincial Natural Science Foundation of China (Grant No. 621RC510). Both M. Fu and H. Li are corresponding authors.} }
\author[a,b]{Jiyun Guo}
\author[a]{Miao Fu \footnote{Corresponding author: miaofu@tju.edu.cn}}
\author[a]{Yuqin Zhang}
\author[b]{Haiyan Li \footnote{Corresponding author: lhy9694@163.com}}
\affil[a]{Center for Applied Mathematics, Tianjin University, Tianjin, 300072, China}
\affil[b]{School of Science, Hainan University, Haikou, 570228, China}
    \maketitle

    \begin{abstract}
The degree sequence of a graph is the sequence of the degrees of its vertices. If $\pi$ is a degree sequence of a graph $G$, then $G$ is a realization of $\pi$ and $G$ realizes $\pi$. Determining when a sequence of positive integers is realizable as a degree sequence of a simple graph has received much attention. One of the early results, by Erd\"{o}s and Gallai, characterized degree sequences of graphs. The result was strengthened by Hakimi and Havel. Another generalization is derived by Cai et al. Hoogeveen and Sierksma listed seven criteria and gave a uniform proof. In addition, Gale and Ryser independently established a characterization by using network flows. We extend Gale and Ryser's results from bipartite graphs to tripartite graphs and even $n$-partite graphs. As corollaries, we give a necessary condition and a sufficient condition for the triple $(\sigma_1, \sigma_2, \sigma_3)$ to be realizable by a tripartite graph, where $\sigma_1$, $\sigma_2$ and $\sigma_3$ are all non-increasing sequences of nonnegative integers. We also give a stronger necessary condition for $(\sigma_1, \sigma_2, \sigma_3)$ to be realizable by a tripartite graph.


        \begin{flushleft}

        \end{flushleft}

    \end{abstract}

    \section{Introduction}

We first introduce some terminology and notations. A sequence $\lambda=(\lambda_1, \cdots, \lambda_n)$ of nonnegative integers is said to be realizable by a graph (i.e., graphic) if it is the degree sequence of a simple graph $G$ on $n$ vertices and such a graph $G$ is called a realization of $\lambda$. The following celebrated theorem due to Erd\"{o}s and Gallai \cite{EG1960} established a characterization of $\lambda$ that is realizable.

 \begin{theorem} (Erd\"{o}s and Gallai \cite{EG1960})\label{11t}
	Let $\lambda=(\lambda_1, \cdots, \lambda_n)$ be a sequence of nonnegative integers with $\lambda_1\ge\cdots\ge \lambda_n$ and $\sum_{i=1}^n\lambda_i\equiv0 (mod~2)$. Then $\lambda$ is realizable by a graph if and only if: for $\mu=0, 1, \cdots, n$, we have
		\begin{equation*}
	\sum_{i=1}^\mu\lambda_i\le  \mu(\mu-1)+\sum_{i=\mu+1}^n\min\{\mu, \lambda_i\}.
	\end{equation*}
\end{theorem}

There have been several proofs of this algebraic characterization, which can be found in \cite{B73,C1986, H1969, T2010}. Cai, Deng and Zang in \cite{C1986}, answering the question of Niessen \cite{N1998}, obtained an extension of the Erd\"{o}s-Gallai result. Yin-Li \cite{YL2005}, Hakimi \cite{H1962} and Havel \cite{H1955}also contributed their own beautiful characterizations. In addition, the notion of realizable sequences (i.e., graphic sequences)  can be extended to bigraphic sequences. There are two sequences of nonnegative integers such that each sequence is the sequence of vertex degrees of a partite set of bipartite graph. Gale \cite{G1957} and Ryser \cite{R1957}, independently, developed a characterization of two sequences that realize a bigraph utilizing results on network flows.

 \begin{theorem} (Gale \cite{G1957}, Ryser \cite{R1957})\label{12t}
	Let $\delta=(\delta_1,\delta_2,\ldots, \delta_m)$ and  $\varphi=(\varphi_1,\varphi_2$, $\ldots, \varphi_n)$ be two sequences of nonnegative integers with $\delta_1\ge \delta_2\ge \cdots\ge \delta_m, \varphi_1\ge \varphi_2\ge \cdots\ge \varphi_n$ and $\sum_{i=1}^m\delta_i=\sum_{i=1}^n\varphi_i$, then  $(\delta; \varphi)$ is realizable by a simple bigraph if and only if

$$\sum_{i=1}^t\varphi_i\le  \sum_{i=1}^m\min\{\delta_i,t\} $$
for $1\le t\le n$.
\end{theorem}

Moreover, the degree sequences for hypergraphs have been most extensively studied, one can refer to \cite{BE2013,DG2021, FPR2021, GL2013}.
And the results of this paper can also be extended to degree sequences of hypergraphs.

\section{Tripartite graphs}

We first consider the existence of tripartite graphs, which can be deduced from Theorem \ref{13t} and Theorem \ref{14t} in the subsequent.

\begin{theorem} \label{13t}
Let $\eta_1=([a_1, b_1], \cdots, [a_m, b_m])$, $\eta_2=([c_1, d_1], \cdots, [c_n, d_n])$ and $\eta_3=([g_1, h_1], \cdots, [g_l, h_l])$ be three sequences of intervals composed of nonnegative integers with $a_1\ge \cdots \ge a_m$, $c_1\ge \cdots \ge c_n$ and  $g_1\ge \cdots \ge g_l$. If they satisfy the following inequalities:
\begin{equation}\label{e1}
\sum_{i=1}^\epsilon\lceil\frac{a_i}{2}\rceil\le  \sum_{j=1}^n\min\{\lfloor\frac{d_j}{2}\rfloor, \epsilon\}  \ \ \mbox{ for every  } \epsilon \mbox{ with } 1\le \epsilon\le m,
\end{equation}
\begin{equation}\label{e2}
\sum_{i=1}^\epsilon\lceil\frac{a_i}{2}\rceil\le  \sum_{j=1}^l\min\{\lfloor\frac{h_j}{2}\rfloor, \epsilon\}  \ \ \mbox{ for every  } \epsilon \mbox{ with } 1\le \epsilon\le m,
\end{equation}	
\begin{equation}\label{e3}
\sum_{i=1}^p\lceil\frac{c_i}{2}\rceil\le  \sum_{j=1}^m\min\{\lfloor\frac{b_j}{2}\rfloor, p\}  \ \ \mbox{ for every  } p \mbox{ with } 1\le p\le n,
\end{equation}
\begin{equation}\label{e4}
\sum_{i=1}^p\lceil\frac{c_i}{2}\rceil\le  \sum_{j=1}^l\min\{\lfloor\frac{h_j}{2}\rfloor, p\}  \ \ \mbox{ for every  } p \mbox{ with } 1\le p\le n,
\end{equation}
 \begin{equation}\label{e5}
\sum_{i=1}^q\lceil\frac{g_i}{2}\rceil\le  \sum_{j=1}^m\min\{\lfloor\frac{b_j}{2}\rfloor, k\}  \ \ \mbox{ for every  } q \mbox{ with } 1\le q\le l,
\end{equation}
\begin{equation}\label{e6}
\sum_{i=1}^q\lceil\frac{g_i}{2}\rceil\le  \sum_{j=1}^n\min\{\lfloor\frac{d_j}{2}\rfloor, k\}  \ \ \mbox{ for every } q \mbox{ with } 1\le q\le l,
\end{equation}
then $(\eta_1; \eta_2; \eta_3)$ is realizable by a tripartite graph, namely, there exists a tripartite graph $G$  with vertex sets $X=\{x_1, \cdots, x_m\}$, $Y=\{y_1, \cdots, y_n\}$, $Z=\{z_1, \cdots, z_l\}$ so that $a_i\le d_G(x_i) \le b_i$ for $i=1, \cdots, m$, $c_i\le d_G(y_i) \le d_i$ for $i=1, \cdots, n$ and $g_i\le d_G(z_i) \le h_i$ for $i=1, \cdots, l$.
\end{theorem}

\pf In \cite{GGT2011}, Garg, Goel and Triparthi developed two results on graphical sequences using a constructive approach. For the proof, we refer to the algorithm in their paper. We consider three sequences of intervals  $\eta_1=([a_1, b_1], \cdots, [a_m, b_m])$, $\eta_2=([c_1, d_1],  \cdots,$  $[c_n, d_n])$ and $\eta_3=([g_1, h_1], \cdots, [g_l, h_l])$ of nonnegative integers with $a_1\ge \cdots \ge a_m$, $c_1\ge \cdots \ge c_n$ and  $g_1\ge \cdots \ge g_l$. Assume  $\eta_1, \eta_2, \eta_3$ satisfy the inequalities (\ref{e1}) to (\ref{e6}). Our goal is to construct a tripartite graph $G$, which will be obtained by constructing three bipartite graphs $G_1$, $G_2$ and $G_3$. We first construct a bipartite graph $G_1^\prime$ with partite sets $\{x_1, \cdots, x_m\}$ and $\{y_1, \cdots, y_n\}$ such that $d_Y(x_i)=\lceil\frac{a_i}{2}\rceil$ for $i=1, \cdots, m$ and $d_X(y_j)\le\lfloor\frac{d_j}{2}\rfloor$ for $j=1, \cdots, n$, where $d_Y(x_i)$ denotes the contribution to the degree of vertex $x_i$ from edges incident to vertices in $Y$. We define the critical index to be the largest index $r$ such that $d_Y(x_i)=\lceil\frac{a_i}{2}\rceil$ for $1\le i < r$ and $d_Y(x_r)<\lceil\frac{a_r}{2}\rceil$. We shall iteratively get rid of the difference $\lceil\frac{a_r}{2}\rceil-d_Y(x_r)$ at vertex $x_r$ while keeping $d_Y(x_i)=\lceil\frac{a_i}{2}\rceil$ for $1\le i <r$ and  $d_X(y_j)\le \lfloor\frac{d_j}{2}\rfloor$ for $1\le j \le n$.  There must be a vertex $v\in N_Y(x_i)\setminus N_Y(x_r)$ for $1\le i< r$, since $d_Y(x_i)=\lceil\frac{a_i}{2}\rceil\ge \lceil\frac{a_r}{2}\rceil>d_Y(x_r)$, where $N_Y(x_i)$ denotes the neighbours set of $x_i$ in $Y$.

	 \begin{enumerate}[fullwidth,labelsep=0.5em,listparindent=2em,label= \textbf{Case \arabic*}]

	\item Suppose, for some $j$, $y_j\leftrightarrow x_k$ for some $k>r$ and $y_j\nleftrightarrow x_\epsilon$ for some $\epsilon\le r$. If $\epsilon=r$, replace $x_ky_j$ by $x_ry_j$. If $\epsilon<r$, replace  $x_ky_j$, $x_\epsilon v$ by $x_\epsilon y_j$, $x_rv$, where $v\in N_Y(x_\epsilon)\setminus N_Y(x_r)$.
	
	\item Suppose for some $j$, $d_X(y_j)<\lfloor\frac{d_j}{2}\rfloor$ and $y_j\nleftrightarrow x_\epsilon$ for some $\epsilon\le r$. If $\epsilon=r$, connect an edge between $x_r$ and $y_j$. If $\epsilon<r$, replace $x_\epsilon v$ by $x_\epsilon y_j$, $x_rv$, where $v\in N_Y(x_\epsilon)\setminus N_Y(x_r)$.
	
	If neither of the above cases will apply, then
	\begin{equation}\label{e7}
	\sum_{i=1}^{r-1}\lceil\frac{a_i}{2}\rceil+d_Y(x_r)= \sum_{i=1}^{r}d_Y(x_i)=\sum_{j=1}^n\min\{d_X(y_j), r\}= \sum_{j=1}^n\min\{\lfloor\frac{d_j}{2}\rfloor, r\}.
	\end{equation}
By (\ref{e1}), (\ref{e7}) and  $d_Y(x_r)\le\lceil\frac{a_r}{2}\rceil$, we have  $d_Y(x_r)=\lceil\frac{a_r}{2}\rceil$. Increasing the value of $r$ by one, and applying the similar procedure yields the bipartite graph $G_1^\prime$, with vertex sets $\{x_1, \cdots, x_m\}$ and $\{y_1, \cdots, y_n\}$ satisfying $d_Y(x_i)=\lceil\frac{a_i}{2}\rceil$ for $i=1, \cdots, m$ and $d_X(y_j)\le\lfloor\frac{d_j}{2}\rfloor$ for $j=1, \cdots, n$.

Next, we shall construct a new bipartite graph $G_1$ on the basis of $G_1^\prime$. Now we define another critical subscript to be the maximal subscript $s$ such that $d_X(y_i)\ge\lceil\frac{c_i}{2}\rceil$ for $1\le i< s$ and $d_X(y_s)<\lceil\frac{c_s}{2}\rceil$. We will get rid of the difference at vertex $y_s$ while keeping  $d_X(y_j)\ge\lceil\frac{c_j}{2}\rceil$ for $1\le j<s$ and $\lceil\frac{a_i}{2}\rceil\le d_Y(x_i)\le \lfloor\frac{b_i}{2}\rfloor$ for $i=1, \cdots, m$ unless there is another situation.

\item If $d_X(y_i)>\lceil\frac{c_i}{2}\rceil$ for $1\le i< s$, then replace $y_iv$ by $y_sv$, where $v\in N_X(y_i)\setminus N_X(y_s)$.

If the above three cases no longer apply, then analogous to Eq (\ref{e7}), we arrive at $d_X(y_s)=\lceil\frac{c_s}{2}\rceil$. Increase $s$ by 1, and applying the same steps yields the required bipartite graph $G_1$ with vertex sets $\{x_1, \cdots, x_m\}$ and $\{y_1, \cdots, y_n\}$ such that $\lceil\frac{a_i}{2}\rceil\le d_Y(x_i)\le \lfloor\frac{b_i}{2}\rfloor$ for $i=1, \cdots, m$ and $d_X(y_i)=\lceil\frac{c_i}{2}\rceil$ for $i=1, \cdots, n$.

To proceed with our proof, we shall construct another bipartite graph $G_2^\prime$ with partite sets $\{x_1, \cdots, x_m\}$ and $\{z_1, \cdots, z_l\}$ satisfying  $d_Z(x_i)=\lceil\frac{a_i}{2}\rceil$ for $1\le i\le m$ and $d_X(z_j)\le\lfloor\frac{h_j}{2}\rfloor$ for $1\le j\le l$. Define the critical subscript to be the largest subscript $r_1$ so that $d_Z(x_i)=\lceil\frac{a_i}{2}\rceil$ for $1\le i< r_1$ and $d_Z(x_{r_1})<\lceil\frac{a_{r_1}}{2}\rceil$. We will repeatedly cut the difference $\lceil\frac{a_{r_1}}{2}\rceil-d_Z(x_{r_1})$ at vertex $x_{r_1}$ while keeping $d_Z(x_{i})=\lceil\frac{a_{i}}{2}\rceil$ for $1\le i< r_1$ and $d_X(z_j)\le\lfloor\frac{h_j}{2}\rfloor$ for $1\le j\le l$. There must be a vertex $v\in N_Z(x_i)\setminus N_Z(x_{r_1})$ for $1\le i < r_1$, since $d_Z(x_{i})=\lceil\frac{a_{i}}{2}\rceil>\lceil\frac{a_{{r_1}}}{2}\rceil>d_Z(x_{r_1})$.

\item  Suppose, for some $j$, $z_j\leftrightarrow x_k$ for some $k>r_1$ and $z_j \nleftrightarrow x_{\epsilon}$ for some $\epsilon\le r_1$. If $\epsilon= r_1$, replace $x_kz_j$ with $x_{r_1}z_j$. If $\epsilon<r_1$, replace $x_kz_j, x_\epsilon v$ with $x_\epsilon z_j, x_{r_1}v$, where $v\in  N_Z(x_\epsilon)\setminus N_Z(x_{r_1})$.

\item Suppose, for some $j$, $d_X(z_j)<\lfloor\frac{h_j}{2}\rfloor$ and $z_j \nleftrightarrow x_{\epsilon}$ for some  $\epsilon\le r_1$. If $\epsilon= r_1$ add the edge $x_{r_1}z_j$. If $\epsilon<r_1$, replace $vx_\epsilon$ by $x_\epsilon z_j, vx_{r_1}$, where $v\in  N_Z(x_\epsilon)\setminus N_Z(x_{r_1})$.

If Case 4 and Case 5 can not apply, then similar to Eq. \ref{e7}, we get $d_Z(x_{r_1})=\lceil\frac{a_{r_1}}{2}\rceil$. Increase $r_1$ by 1 and continue the steps, then we can obtain the bigraph $G_2^\prime$ with partite sets $\{x_1, \cdots, x_m\}$ and $\{z_1, \cdots, z_l\}$ satisfying $d_Z(x_{i})=\lceil\frac{a_{i}}{2}\rceil$ for $i=1, \cdots, m$ and $d_X(z_{j})\le\lfloor\frac{h_{j}}{2}\rfloor$ for $j=1, \cdots, l$.

Now let $s_1$ be the largest index so that $d_X(z_j)\ge\lceil\frac{g_j}{2}\rceil$ for $1\le j<s_1$ and $d_X(z_{s_1})<\lceil\frac{g_{s_1}}{2}\rceil$. We decrease the deficiency at vertex $z_{s_1}$ while maintaining $d_X(z_j)\ge\lceil\frac{g_j}{2}\rceil$ for $1\le j<s_1$, and  $\lceil\frac{a_i}{2}\rceil\le d_Z(x_i)\le \lfloor\frac{b_i}{2}\rfloor$ for $1\le i\le m$, except that an additional case arises.

\item If $d_X(z_i)>\lceil\frac{g_i}{2}\rceil$ for $1\le i<s_1$, then replace $vz_i$ by $vz_{s_1}$, where $v\in N_X(z_i)\setminus N_X(z_{s_1})$.

If Case 4, Case 5 and Case 6 can not apply, then analogous to Eq. \ref{e7}, we have $d_X(z_{s_1})=\lceil\frac{g_{s_1}}{2}\rceil$. Increase $s_1$ by 1, and continue the steps, then we obtain the bigraph $G_2$ with partite set $\{x_1, \cdots, x_m\}$ and $\{z_1, \cdots, z_l\}$ satisfying that $\lceil\frac{a_i}{2}\rceil\le d_Z(x_i)\le \lfloor\frac{b_i}{2}\rfloor$ for $i=1, 2, \cdots, m$ and $d_X(z_i)=\lceil\frac{g_i}{2}\rceil$ for $i=1, \cdots, l$. Next, we shall construct a new bigraph $G_3^\prime$ with partite sets  $\{y_1, \cdots, y_n\}$ and $\{z_1, \cdots, z_l\}$ satisfying $d_Z(y_i)=\lceil\frac{c_i}{2}\rceil$ for $1\le i\le n$ and $d_Y(z_j)\le\lfloor\frac{h_j}{2}\rfloor$ for $j=1, \cdots, l$. We define the maximal critical subscript  $r_2$ such that $d_Z(y_i)=\lceil\frac{c_i}{2}\rceil$ for $1\le i< r_2$ and $d_Z(y_{r_2})<\lceil\frac{c_{r_2}}{2}\rceil$. We shall repeatedly remove the difference $\lceil\frac{c_{r_2}}{2}\rceil-d_Z(y_{r_2})$ at vertex $y_{r_2}$ while keeping $d_Z(y_i)=\lceil\frac{c_i}{2}\rceil$ for $1\le i< r_2$ and $d_Y(z_j)\le\lfloor\frac{h_j}{2}\rfloor$ for $1\le j\le l$. Notice that there exists a vertex $v\in N_Z(y_i)\setminus N_Z(y_{r_2})$ for $1\le i\le r_2$, since $d_Z(y_i)=\lceil\frac{c_i}{2}\rceil\ge \lceil\frac{c_{r_2}}{2}\rceil>d_Z(y_{r_2})$.

\item Suppose for some $j$, $z_j\leftrightarrow y_k$ for some $k>r_2$ and $z_j\nleftrightarrow y_\epsilon$ for some $\epsilon\le r_2$. If $\epsilon= r_2$, replace $y_kz_j$ by $z_jy_{r_2}$. If $\epsilon<r_2$, replace $y_kz_j$, $vy_\epsilon$ by $y_\epsilon z_j$, $vy_{r_2}$, where $v\in N_Z(y_\epsilon)\setminus N_Z(y_{r_2})$.

\item  Suppose for some $j$, $d_Y(z_j)<\lfloor\frac{h_j}{2}\rfloor$ and $z_j\nleftrightarrow y_\epsilon$ for some $\epsilon\le r_2$. If $\epsilon= r_2$, add the edge $y_{r_2}z_j$. If $\epsilon<r_2$, replace $vy_\epsilon$ by $y_\epsilon z_j$, $vy_{r_2}$, where $v\in N_Z(y_\epsilon)\setminus N_Z(y_{r_2})$.

If Case 7 and Case 8 can not apply, then
\begin{equation*}\sum_{i=1}^{r_2-1}\lceil\frac{c_i}{2}\rceil+d_Z(y_{r_2})=\sum_{i=1}^{r_2}d_Z(y_i)=\sum_{j=1}^l\min\{d_Y(z_j), r_2\}=\sum_{j=1}^l\min\{\lfloor\frac{h_j}{2}\rfloor, r_2\}.
\end{equation*}
By $(5)$ and $d_Z(y_i)=\lceil\frac{c_i}{2}\rceil$ for $1\le i<r_2$, we have $d_Z(y_{r_2})=\lceil\frac{c_{r_2}}{2}\rceil$. Increasing $r_2$ by 1, and applying the same steps reduces to a bipartite graph $G_3^\prime$, with partite sets $\{y_1,\cdots, y_n\}$ and $\{z_1, \cdots, z_l\}$ satisfying $d_Z(y_i)=\lceil\frac{c_i}{2}\rceil$ for $1\le i\le n$ and $d_Y(z_j)\le\lfloor\frac{h_j}{2}\rfloor$ for $1\le j\le l$.

We define the maximal critical subscript $t$ such that $d_Y(z_i)\ge\lfloor\frac{g_i}{2}\rfloor$ for $1\le i< t$ and $d_Y(z_t)<\lceil\frac{g_t}{2}\rceil$. We will get rid of the difference at vertex $z_t$ while keeping $\lceil\frac{c_i}{2}\rceil\le d_Z(y_i)\le \lfloor\frac{d_i}{2}\rfloor$ for $1\le i\le n$ and $d_Y(z_j)\ge\lfloor\frac{g_j}{2}\rfloor$ for $1\le j<t$ unless there is a new case.

\item If $d_Y(z_j)>\lceil\frac{g_j}{2}\rceil$ for some $j<t$, then replace $vz_j$ by $vz_t$, where $v\in N_Y(z_j)\setminus N_Y(z_t)$.

If Case 7-9 can't apply, then analogous to Eq. \ref{e7}, we get that $d_Y(z_t)=\lceil\frac{g_t}{2}\rceil$. Increase $t$ by one, and continue. The argument above gives rise to a bipartite graph $G_3$ with vertex sets $\{y_1,\cdots, y_n\}$ and $\{z_1, \cdots, z_l\}$ with $\lceil\frac{c_i}{2}\rceil\le d_Z(y_i)\le \lfloor\frac{d_i}{2}\rfloor$ for $i=1, \cdots, n$ and $d_Y(z_j)=\lceil\frac{g_j}{2}\rceil$ for $j=1, \cdots, l$.

From what has been discussed, we get three bipartite graphs called $G_1, G_2$ and $G_3$, where

\begin{equation*}G_1: \lceil\frac{a_i}{2}\rceil\le d_Y(x_i)\le \lfloor\frac{b_i}{2}\rfloor \mbox{ for } 1\le i\le m, \\  d_X(y_i)=\lceil\frac{c_i}{2}\rceil \mbox{ for } 1\le i\le n;
\end{equation*}

\begin{equation*}G_2: \lceil\frac{a_i}{2}\rceil\le d_Z(x_i)\le \lfloor\frac{b_i}{2}\rfloor \mbox{ for } 1\le i\le m, \\  d_X(z_j)=\lceil\frac{g_j}{2}\rceil \mbox{ for } 1\le j\le l;
\end{equation*}

\begin{equation*}G_3: \lceil\frac{c_i}{2}\rceil\le d_Z(y_i)\le \lfloor\frac{d_i}{2}\rfloor \mbox{ for } 1\le i\le n, \\  d_Y(z_j)=\lceil\frac{g_j}{2}\rceil \mbox{ for } 1\le j\le l.
\end{equation*}

Set $G=G_1\cup G_2 \cup G_3$, then $G$ is the required tripartite graph, since $d_G(x_i)=d_Y(x_i)+d_Z(x_i)$, $d_G(y_i)=d_X(y_i)+d_Z(y_i)$ and $d_G(z_i)=d_X(z_i)+d_Y(z_i)$. Clearly, $2\lceil\frac{a_i}{2}\rceil\le d_G(x_i)\le 2\lfloor\frac{b_i}{2}\rfloor$ for $1\le i\le m$, $2\lceil\frac{c_i}{2}\rceil\le d_G(y_i)\le \lceil\frac{c_i}{2}\rceil+\lfloor\frac{d_i}{2}\rfloor$ for $1\le i\le n$, and $d_G(z_i)=2\lceil\frac{g_i}{2}\rceil\in [g_i, h_i]$ for $1\le i\le l$. One can see $2\lceil\frac{a}{2}\rceil\ge a\ge 2\lfloor\frac{a}{2}\rfloor$ for any positive integer $a$. Thus $G$ is the graph that satisfies the conditions and so we are done.
\end{enumerate}	 \qed

The conditions (\ref{e1})-(\ref{e6}) of Theorem \ref*{13t} are not necessary, as can be seen by taking $\eta_1=([2,3], [0, 2])$, $\eta_2=([2,4], [1, 2])$ and $\eta_3=([1,2], [0, 1])$, which satisfy $a_1\ge a_2$, $c_1\ge c_2$ and $g_1\ge g_2$. It is easy to check that the tripartite graph below is a realization of $(\eta_1; \eta_2; \eta_3)$.
However, (\ref{e5}) does not hold for $k=2$.
Theorem \ref{14t} states the necessary condition for $(\eta_1; \eta_2; \eta_3)$ to be realizable by a tripartite graph.
\begin{figure}[htbp]
	\begin{center}
		\includegraphics[scale=0.7]{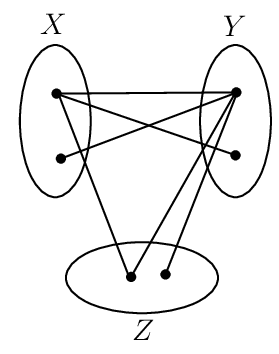}
	\label{Fig.1}
	\end{center}
\end{figure}
  \begin{theorem} \label{14t}
	Let $\eta_1=([a_1, b_1], \cdots, [a_m, b_m])$, $\eta_2=([c_1, d_1], \cdots, [c_n, d_n])$ and $\eta_3=([g_1, h_1], \cdots, [g_l, h_l])$ be three lists of intervals composed of nonnegative integers with $a_1\ge \cdots \ge a_m$, $c_1\ge \cdots \ge c_n$ and  $g_1\ge \cdots \ge g_l$. If $(\eta_1; \eta_2; \eta_3)$ is realizable by a tripartite graph, then
	\begin{equation} \label{e8}
	\sum_{i=1}^\tau a_i\le  \sum_{j=1}^n\min\{d_j, \tau\}+\sum_{j=1}^l\min\{h_j, \tau\}  \ \ \mbox{ for every  } \tau \mbox{ with } 1\le \tau\le m;
	\end{equation} \label{e9}
	\begin{equation}
	\sum_{i=1}^\tau c_i\le  \sum_{j=1}^n\min\{b_j, \tau\}+\sum_{j=1}^l\min\{h_j, \tau\}  \ \ \mbox{ for every  } \tau \mbox{ with } 1\le \tau\le n;
	\end{equation}	
	\begin{equation}\label{e10}
	\sum_{i=1}^\tau g_i\le  \sum_{j=1}^n\min\{b_j, \tau\}+\sum_{j=1}^l\min\{d_j, \tau\}  \ \ \mbox{ for every  } \tau \mbox{ with } 1\le \tau\le l.
	\end{equation}

\end{theorem}

\pf Assume that $G$ is a tripartite graph with partite sets $X, Y$ and $Z$ satisfying the given conditions. Consider the edges incident to a set of $\tau$ vertices in $X$. Each $y_j\in Y$ is incident not more than $\tau$ of these vertices, and also incident to not more than $d_{Y_X}(y_j)$ of these vertices, where $d_{Y_X}(y_j)$ is the maximum contribution to these $\tau$ vertices in $X$ from the edges incident to $y_j\in Y$. Analogous to the case of vertices in $Z$. Thus,
 \begin{eqnarray*}
\sum_{i=1}^\tau a_i\le \sum_{i=1}^\tau d_G(x_i) &\le& \sum_{j=1}^n\min\{d_{Y_X}(y_j), \tau\}+\sum_{j=1}^l\min\{d_{Z_X}(z_j), \tau\}\\
	&\le& \sum_{j=1}^n\min\{d_{G}(y_j), \tau\}+\sum_{j=1}^l\min\{d_{G}(z_j), \tau\}\\
	&\le& \sum_{j=1}^n\min\{d_j, \tau\}+\sum_{j=1}^l\min\{h_j, \tau\}.
\end{eqnarray*}
So, (\ref{e8}) follows.

The other two inequalities can be shown to be true in the same way as desired.\qed

Theorem \ref{13t} and Theorem \ref{14t} reduce to two corollaries, one of which is sufficient and the other one necessary, for $(A; B; C)$ to be realizable by a tripartite graph when $a_i=b_i$ for all $i$, $c_j=d_j$ for all $j$ and $g_k=h_k$ for all $k$, where $A=(b_1, \cdots, b_m)$, $B=(d_1, \cdots, d_n)$ and $C=(h_1, \cdots, h_l)$. Furthermore, the two corollaries extend Gale-Ryser Theorem from a bipartite graph to a tripartite graph.

 \begin{corollary} \label{23c}
	Let $\sigma_1=(a_1, \cdots, a_m)$, $\sigma_2=(b_1, \cdots, b_n)$ and  $\sigma_3=(c_1, \cdots, c_l)$ be three non-increasing sequences of nonnegative integers. If
	
	\begin{equation*}
	\sum_{i=1}^s \lceil {a_i\over 2}\rceil \le  \sum_{j=1}^n\min\{ \lfloor{b_j\over 2}\rfloor, s\}  \ \ \mbox{ for every  } s \mbox{ with } 1\le s\le m
	\end{equation*}
		and
	\begin{equation*}
\sum_{i=1}^t \lceil {a_i\over 2}\rceil \le  \sum_{j=1}^l\min\{ \lfloor{c_j\over 2}\rfloor, t\}  \ \ \mbox{ for every } t \mbox{ with } 1\le t\le m;
\end{equation*}	
all hold, then there exists a tripartite graph $G$ with tripartition $\{(x_1, \cdots, x_m), (y_1, \cdots, y_n), (z_1, \cdots, z_l)\}$ such that $d_G(x_i)=a_i$ for $i=1, \cdots, m$, $d_G(y_j)=b_j$ for $j=1, \cdots, n$ and $d_G(z_k)=c_k$ for $k=1, \cdots, l$, i.e., $(\sigma_1; \sigma_2; \sigma_3)$ is realizable by a tripartite graph $G$.
\end{corollary}

  \begin{corollary} \label{24c}
 	Let $\sigma_1=(a_1, \cdots, a_m)$, $\sigma_2=(b_1, \cdots, b_n)$ and  $\sigma_3=(c_1, \cdots, c_l)$ be three non-increasing sequences of nonnegative integers. If $(\sigma_1; \sigma_2; \sigma_3)$ is realizable by a tripartite graph, then for each $s$ with  $1\le s \le m$,
 	
 	\begin{equation*}
 		\sum_{i=1}^s {a_i} \le  \sum_{j=1}^n\min\{ b_j, s\} +  \sum_{k=1}^l\min\{c_k, s\}.
 	\end{equation*}	
  \end{corollary}

The following theorem gives a necessary condition for the triple $(\sigma_1; \sigma_2; \sigma_3)$ to be realizable by a tripartite graph, and the condition is stronger than the one in Corollary \ref{24c}.

  \begin{theorem} \label{25t}
	Let $\sigma_1=(a_1, \cdots, a_m)$, $\sigma_2=(b_1, \cdots, b_n)$ and  $\sigma_3=(c_1, \cdots, c_l)$ be three non-increasing sequences of nonnegative integers. If $(\sigma_1; \sigma_2; \sigma_3)$ is realizable by a tripartite graph, then for each $\delta$ with  $1\le \delta \le m$,
	
	\begin{equation*}
	\sum_{i=1}^\delta  {a_i} \le  \min\{\sum_{j=1}^n b_j-\mu, n\delta\}+\min\{\sum_{k=1}^l c_k-\mu, l\delta\},
	\end{equation*}	
where $\mu={1\over2}(\sum_{j=1}^n {b_j}+\sum_{k=1}^l {c_k}-\sum_{i=1}^m {a_i})$.
\end{theorem}

\pf
Suppose that $(\sigma_1; \sigma_2; \sigma_3)$ is realizable by a tripartite graph $G$ with partite sets $X=\{x_1, \cdots, x_m\}$, $Y=\{y_1, \cdots, y_n\}$ and $Z=\{z_1, \cdots, z_l\}$, then $d_G(x_i)=a_i$ for $i=1, \cdots, m$, $d_G(y_j)=b_j$ for $j=1, \cdots, n$ and $d_G(z_k)=c_k$ for $k=1, \cdots, l$. Denote $\mu$ by the number of edges of $G$ between $Y$ and $Z$. Then we see that
 \begin{eqnarray*}
	\mu &=& \sum_{\epsilon=1}^nd(y_\epsilon)-\sum_{\epsilon=1}^nd_X(y_\epsilon)\\
	&=& \sum_{\epsilon=1}^nb_\epsilon-\sum_{\epsilon=1}^nd_X(y_\epsilon),
\end{eqnarray*}
and
 \begin{eqnarray*}
	\mu &=& \sum_{\tau=1}^ld(z_\tau)-\sum_{\tau=1}^ld_X(z_\tau)\\
	&=& \sum_{\tau=1}^lc_\tau-\sum_{\tau=1}^ld_X(z_\tau),
\end{eqnarray*}
where $d_X(y_j)$ denotes the contribution of all vertices in $X$ to the vertex $y_j\in Y$. Thus
	\begin{equation*}
\sum_{j=1}^n {d_X(y_j)} = \sum_{j=1}^n  b_j -  \sum_{k=1}^l c_k+\sum_{k=1}^l {d_X(z_k)},
\end{equation*}	
and so
 \begin{eqnarray*}
	\sum_{i=1}^m a_i &=& \sum_{i=1}^md(x_i)\\ &=&\sum_{j=1}^nd_X(y_j)+\sum_{k=1}^ld_X(z_k)\\
	&=& \sum_{j=1}^nb_j-\sum_{k=1}^lc_k+2\sum_{k=1}^ld_X(z_k).
\end{eqnarray*}

After the transpose, we can obtain
	\begin{equation*}
\sum_{k=1}^l {d_X(z_k)} ={1\over 2} \big(\sum_{i=1}^m {a_i}-\sum_{j=1}^n  b_j+ \sum_{k=1}^l c_k\big),
\end{equation*}
and then
	\begin{equation*}
\mu ={1\over 2} \big(\sum_{j=1}^n  b_j+ \sum_{k=1}^l c_k-\sum_{i=1}^m {a_i}\big).
\end{equation*}

Thus,	\begin{equation*}
\sum_{j=1}^n d_X(y_j)=\sum_{j=1}^n b_j-\mu, \sum_{k=1}^l d_X(z_k)=\sum_{k=1}^l c_k-\mu.
\end{equation*}

Hence, for each $\delta$ with $1\le \delta \le m$, we have
 \begin{eqnarray*}
	\sum_{i=1}^\delta a_i &=& \sum_{i=1}^\delta d(x_i)\\ &\le&\sum_{j=1}^n\min\{d_X(y_j), \delta\}+\sum_{k=1}^l\min\{d_X(z_k), \delta\}\\
	&\le&\min\{\sum_{j=1}^n d_X(y_i), n\delta\}+\min\{\sum_{k=1}^l d_X(z_k), l\delta\}\\
&=& \min\{\sum_{j=1}^n b_j-\mu, n\delta\}+\min\{\sum_{k=1}^l c_k-\mu, l\delta\},
\end{eqnarray*}
where the expression for $\mu$ is shown above.\qed

Based on the previous results, we want to establish a necessary and sufficient condition for $(\sigma_1, \sigma_2, \sigma_3)$ to be realizable by a tripartite graph or a hypergraph. However, it seems difficult to give such a characterization.

 \begin{problem} \label{26p} Determine a characterization of tripartite graphic sequences.
\end{problem}

 \begin{problem} \label{27p} Investigate necessary condition and sufficient condition for $(\sigma_1, \sigma_2, \sigma_3)$ to be realizable by an $r$-uniform hypergraph.
 \end{problem}

\section{$n$-partite graphs}

Motivated by Theorem  \ref{13t} and Theorem \ref{14t}, we also consider the situation of $n$-partite graph, namely, Theorem \ref{15t} and Theorem \ref{16t}.

  \begin{theorem} \label{15t}
	Let $L_1=([a_1^1, b_1^1], \cdots, [a_{k_1}^1, b_{k_1}^1])$, $L_2=([a_1^2, b_1^2], \cdots, [a_{k_2}^2, b_{k_2}^2])$, $\cdots$, $L_n=([a_1^n, b_1^n], \cdots, [a_{k_n}^n, b_{k_n}^n])$ be $n$ sequences of intervals composed of nonnegative integers with $a_1^1\ge a_2^1\ge \cdots \ge a_{k_1}^1$, $\cdots\cdots$, $a_1^n\ge a_2^n\ge \cdots \ge a_{k_n}^n$. If the following $\frac{n(n-1)}{2}$ inequalities hold:

	\begin{equation*}
	\sum_{i=1}^{r_1} \big\lceil \frac{a_i^1}{n-1}\big\rceil\le \min\bigg\{ \sum_{j=1}^{k_s}\min\big\{\big\lfloor \frac{a_j^s}{n-1}\big\rfloor, r_1\big\}, s=2, 3, \cdots, n\bigg\}, 1\le r_1\le k_1;
	\end{equation*}
	$$\cdots\cdots$$	
	\begin{equation*}
\sum_{i=1}^{r_m} \big\lceil \frac{a_i^m}{n-1}\big\rceil\le \min\bigg\{ \sum_{j=1}^{k_s}\min\big\{\big\lfloor \frac{a_j^s}{n-1}\big\rfloor, r_m\big\}, s=1, \cdots, m-1, m+1, \cdots,  n\bigg\}, 1\le r_m\le k_m;
\end{equation*}
		$$\cdots\cdots$$	
	\begin{equation*}
	\sum_{i=1}^{r_n} \big\lceil \frac{a_i^n}{n-1}\big\rceil\le \min\bigg\{ \sum_{j=1}^{k_s}\min\big\{\big\lfloor \frac{a_j^s}{n-1}\big\rfloor, r_n\big\}, s=1, \cdots,  n-1\bigg\}, 1\le r_n\le k_n,
	\end{equation*}
	then $(L_1; \cdots; L_n)$ is realizable  by a $n$-partite graph.
\end{theorem}

\begin{theorem} \label{16t}
If $(L_1; \cdots; L_n)$ is realizable by a $n$-partite graph, then
	
	\begin{equation*}
\sum_{i=1}^w a_i^1\le  \sum_{j=1}^{k_2}\min\{b_j^2, w\}+\cdots+\sum_{j=1}^{k_n}\min\{b_j^n, w\}, 1\le w\le k_1;
\end{equation*}
	$$\cdots\cdots$$	
	\begin{equation*}
\sum_{i=1}^w a_i^n\le  \sum_{j=1}^{k_1}\min\{b_j^1, w\}+\cdots+\sum_{j=1}^{k_{n-1}}\min\{b_j^{n-1}, w\}, 1\le w\le k_n.
\end{equation*}
\end{theorem}

\noindent\textbf{ Proof of Theorem \ref{15t}} By induction on $n$, when $n=2$ and $n=3$, Theorem \ref{15t} is trivially true. Now suppose that it holds for all $n\le t$, then by the induction hypothesis, we have a $t$-partite graph $G$ satisfying the inequalities in Theorem \ref{15t}. Similar to the construction of tripartite graph, we can construct a $t+1$-partite graph based on $G$.\qed

\bigskip

\noindent\textbf{ Proof of Theorem \ref{16t}} By symmetry, we just have to prove the first inequality by induction on $n$. When $n=2$, it is clearly true. Now suppose that the inequality is true for all $n$ values less that $t$. Thus, for each $w$ with $1\le w\le k_1$, we have
 \begin{eqnarray*}
\sum_{i=1}^w a_i^1 &\le&  \sum_{j=1}^{k_2}\min\{b_j^2, w\}+\cdots+\sum_{j=1}^{k_t}\min\{b_j^t, w\}\\ &<& \sum_{j=1}^{k_2}\min\{b_j^2, w\}+\cdots+\sum_{j=1}^{k_t}\min\{b_j^t, w\}+\sum_{j=1}^{k_{t+1}}\min\{b_j^{t+1}, w\}.
\end{eqnarray*}
So when $n=t+1$, the inequality is also true, and hence the proof is completed.
 \qed

\begin{corollary} \label{33c}
Let $\sigma_1=(a_1^1, a_2^1, \cdots, a_{k_1}^1)$, $\sigma_2=(a_1^2, a_2^2, \cdots, a_{k_2}^2)$, $\cdots$, $\sigma_n=(a_1^n, a_2^n, \cdots, a_{k_n}^n)$ be $n$ non-increasing sequences of nonnegative integers. If for each $r$ with $1\le r \le k_1$,
	
	\begin{equation*}
	\sum_{i=1}^r \big\lceil {a_i^1\over n-1}\big\rceil\le  \min\big\{\sum_{i=1}^{k_s}\min\{\big\lfloor {a_j^s\over n-1}\big\rfloor,r\}, s=2, 3, \cdots, n\big\},
	\end{equation*}
then there exists a $n$-partite graph $G$ with partition $X_1=\{x_1^1, \cdots, x_{k_1}^1\}$, $\cdots$,  $X_n=\{x_1^n, \cdots, x_{k_n}^n\}$ so that their degree sequences are $\sigma_1, \sigma_2, \cdots, \sigma_n$ respectively.
\end{corollary}

\begin{corollary} \label{34c}
	Let $\sigma_1=(a_1^1, a_2^1, \cdots, a_{k_1}^1)$,  $\cdots$, $\sigma_n=(a_1^n, a_2^n, \cdots, a_{k_n}^n)$ be $n$ non-increasing sequences of nonnegative integers. If $(\sigma_1, \sigma_2, \cdots, \sigma_n)$ is realizable by a $n$-partite graph, then for each $w$ with $1\le w \le k_1$,
	
	\begin{equation*}
	\sum_{i=1}^w a_i^1\le \sum_{j=1}^{k_2} \min\{a_j^2, w\}+\cdots+\sum_{j=1}^{k_n} \min\{a_j^n, w\} .
	\end{equation*}

\end{corollary}

\end{document}